\documentclass[12pt]{amsart}

\input xypic.tex
\usepackage{verbatim,amsmath,latexsym,amssymb}
\newcounter{tempsection}

\theoremstyle{plain}
\newtheorem{thm}{Theorem}[section]
\newtheorem{theorem}{Theorem}[tempsection]
\newtheorem{lem}[thm]{Lemma}
\newtheorem{prop}[thm]{Proposition}
\newtheorem{cor}[thm]{Corollary}

\theoremstyle{definition}
\newtheorem{alg}[thm]{Algorithm}
\newtheorem{df}[thm]{Definition}

\newtheorem{rem}[thm]{Remark}

\def\cd{\operatorname{cd}}

\def\gr{\operatorname{gr}}

\def\im{\operatorname{im}}

\def\sgn{\operatorname{sgn}}

\def\tor{\operatorname{Tor}}
\def\var{\operatorname{Var}}

\def\m{{\mathfrak m}}
\def\qed{\hfill$\Box$}
\def\del{\partial}

\def\C{{\mathbb C}}	
\def\D{{\mathcal D}}
\def\E{{\mathcal E}}
\def\F{{\mathcal F}}
\def\I{{\mathcal I}}
\def\J{{\mathcal J}}
\def\K{{\mathcal K}}
\def\L{{\mathcal L}}
\def\M{{\mathcal M}}

\def\MV{{\mathcal M\mathcal V}}
\def\O{{\mathcal O}}	
\def\Z{{\mathbb Z}}
	
\def\ilim{{\lim\limits_{\longleftarrow}}}
\def\into{\hookrightarrow}
\def\onto{\to\hskip-1.7ex\to}

\def\mylabel#1{\label{#1}}
\def\ignore#1{}

\numberwithin{equation}{section}
\begin{document}

\title[Computation of de Rham Cohomology]
       {Algorithmic Computation of de Rham Cohomology of Complements
       of Complex Affine Varieties}
\author[Uli Walther]{Uli Walther\\ University of Minnesota}
\address{School of Mathematics, University of Minnesota, Minneapolis, MN 55455}
\email{walther@math.umn.edu}
\begin{abstract}
Let $X=\C^n$. In this paper we present an algorithm that computes the
de Rham 
cohomology groups $H^i_{dR}(U,\C)$ where $U$ is the complement of an
arbitrary Zariski-closed set $Y$ in $X$. 

Our
algorithm is a merger of the algorithm given by T.~Oaku and
N.~Takayama (\cite{O-T2}), 
who considered the case where $Y$ is a hypersurface, and
our methods from \cite{W-1} for the computation of local cohomology. 
We further extend the algorithm to compute de Rham cohomology groups
with support $H^i_{dR,Z}(U,\C)$ where again $U$ is an arbitrary
Zariski-open subset of $X$ and $Z$ is an arbitrary Zariski-closed
subset of $U$. 

Our main tool is the generalization of the restriction process from
\cite{O-T1} to complexes of modules over the Weyl algebra.

All presented algorithms are 
based on Gr\"obner basis computations in the Weyl algebra.
\end{abstract}

\maketitle

\tableofcontents

\section{Introduction} 
In his famous paper \cite{DRCAV} R.\ Hartshorne introduced the concept
of algebraic de Rham cohomology of algebraic varieties 
as an analog to classical (singular)
cohomology and proved, using results of A.~Grothendieck and
P.~Deligne, that it agrees with classical cohomology if the base field
is $\C$. Moreover, he also defined the notion of algebraic de Rham
cohomology with supports and proved that it fits into certain natural
long exact 
sequences related to inclusion maps (\ref{long-local-seq}). 

In \cite{O-T2}, the authors give an algorithm that computes (by
Gr\"obner basis computations) the algebraic de Rham cohomology of the
complement $U$ of any given  hypersurface $Y$ of $X=\C^n$. Their method is 
based on the initial definition of Hartshorne, as the hypercohomology of
the de Rham complex on  
$U$. They show that this complex is in the
derived category the same as the tensor product over $\O_X$ of the sheaf of
differential $n$-forms on $X$ with a resolution of $\O_U$, $\O_U$
considered  as a module over
the sheaf of differential operators on $X$. The computation of the
hypercohomology of the latter complex reduces to computation of usual
cohomology of the global sections since $U$ is affine and the sheaves
involved are quasi-coherent. An algorithm to compute the cohomology
of complexes of the type one gets after taking global sections was
 given in \cite{O-T1}. The strategy is use the method of restriction
of a $D$-module to a linear subvariety (\cite{O-T1}, section 5). 

In this note we shall prove
\setcounter{tempsection}{4}
\setcounter{theorem}{2}
\begin{theorem}
The de Rham cohomology groups of the complement of an affine complex
variety are effectively computable by means of Gr\"obner basis
computations in rings of differential operators.
\end{theorem}

In fact, we shall 
first generalize the restriction process to the
restriction of a 
complex to a linear subvariety. Then, as applications,  we obtain an algorithm
that computes de Rham cohomology of arbitrary Zariski-open 
$U$, and 
an algorithm that computes 
de Rham cohomology of Zariski open sets 
with support in a Zariski closed subset $Z$ of $U$.

Now we shall give an overview of the structure of this paper. Let
$D_n=\C[x_1,\ldots,x_n]\langle \del_1,\ldots,\del_n\rangle$, the $n$th
Weyl algebra over $\C$.

First of all, in section \ref{section-global}, we show that
if $U$ is the complement of any Zariski closed set $Y$ defined by
$f_1,\ldots,f_r$ in $X$ then 
computation of the de Rham cohomology of $U$ can be performed by
computing the cohomology of the tensor product over $D_n$ of a $D_n$-free
resolution of $D_n/(\del_1,\ldots,\del_n)D_n$ with the Mayer-Vietoris
complex $MV^\bullet(F)$ associated to $f_1,\ldots,f_r$ (cf.\
subsection \ref{MV}).    

In the following section we compute a certain
$D_n$-free complex that is quasi-isomorphic to $MV^\bullet(F)$. In fact,
we present a method that computes for an arbitrary complex of
finitely generated $D_n$-modules $C^\bullet$ with cohomology  
that is specializable to the origin a $D_n$-free complex $A^\bullet$ that is
quasi-isomorphic to $C^\bullet$  and has certain properties related to the
$V_n$-filtration (for facts about the $V_n$-filtration, see also \cite{O-T1}).

Section \ref{restrict-complex} is devoted to the explicit computation of the
derived tensor product $D_n/(\del_1,\ldots,\del_n)\cdot D_n
\otimes^L_{D_n}C^\bullet$
where again $C^\bullet$ is required to have specializable
cohomology, but otherwise is arbitrary. 
A corollary of this computation will be an algorithm that
computes $H^\bullet_{dR}(U,\C)$, (\ref{alg-de-rham}).

In section \ref{section-local} we review the definition of algebraic
de Rham cohomology with supports and exhibit an algorithm that computes
$H^i_{dR,Z}(X\backslash Y,\C)$ for arbitrary subvarieties $Y,Z$ of
$X$. The idea here is similar to the original argument, twisted with
the \v Cech complex associated to $Z$. 


\section{Algebraic de Rham Cohomology}
\mylabel{section-global}
\subsection{Notation} Throughout this article, we shall use the following
notation. $\C$ will stand for the field of complex numbers, $X$
denotes affine $n$-dimensional $\C$-space $\C^n$ and $Y$ will be a
subvariety of $X$ cut out by polynomials $\{f_1,\ldots,f_r\}\subseteq
R$ where $R$ is $\C[x_1,\ldots,x_n]$.  
Let $U=X\backslash Y$.

$D_n$ will be the ring of differential operators on $X$ (also called the
Weyl algebra) generated by the multiplications by the $x_i$ (which we
will call also $x_i$) and the partial derivatives $\del_i=\del/\del
x_i$. Set $\O_X$ to be the structure sheaf on $X$. 
${\mathcal D}_X$ will be the sheaf version of $D_n$, 
$\D_X=\O_X\otimes_RD_n$. We set $\Omega=D_n/(\del_1,\ldots,\del_n)D_n$ and
$\Omega(\D)=\Omega\otimes_{D_n}\D$. 

If $\M$ is a $\D$- or $D_n$-module, $\Omega^\bullet(\M)$ will throughout stand
for the de Rham complex of 
$\M$. In other words, $\Omega^k(\M)=\M\otimes_\Z\bigwedge^k(\Z^n)$
and the differential $d$ is defined in the usual way:
$d(u\otimes dx_{i_1}\wedge\ldots\wedge
dx_{i_k})=\sum_{j=1}^n\del_{j}\cdot u\otimes dx_j\wedge 
dx_{i_1}\wedge\ldots\wedge
dx_{i_k}$. 
If $V$ is a variety, $\Omega^\bullet(V)$ will denote the de Rham
complex on $V$.

Furthermore, set $\Omega^\bullet=\Gamma(X,\Omega^\bullet(X))$.

\subsection{Definition of de Rham cohomology}
Recall the idea of completion $\hat \F$ of a quasi-coherent sheaf
$\F$ on $A$ with respect to the closed subset $B$:
if $V$ is open in $B$, 
$\hat\F(V)=\ilim_k(\F(V)/\I^k(V)\cdot\F(V))$ where $\I$ is the sheaf
of ideals defining $B$. 

Algebraic de Rham cohomology of an arbitrary closed subset $B$ of an
arbitrary smooth scheme $A$ over any field $K$ is defined as the
hypercohomology of the complex $\hat\Omega^\bullet(A)$
where the hat denotes completion of $\Omega^\bullet(A)$ with respect
to the system of ideals 
${\mathcal I}$ which defines $B$ in $A$. (For a precise definition of the
maps in $\hat 
\Omega^\bullet(A)$ see \cite{DRCAV}, page 22.)
It is shown in \cite{DRCAV} that
this definition does not depend on the embedding of $B$ in $A$ nor in
fact on $A$ itself.

In the special case where $B$ is smooth, we may take $B=A$ and then
the sheaf of ideals ${\mathcal I}$ is the zero sheaf. In particular, for
open subsets of $X$, $H^i_{dR}(U,\C)$ is the hypercohomology of the
complex $\Omega^\bullet(U)$. 

\subsection{The Idea of Oaku and Takayama}
For this subsection, assume that $Y$ is a hypersurface.

The basic observation is the following
\begin{lem}
The complex  $\Omega^\bullet(D_n)=\Omega^\bullet$ 
is (a complex in the category of
right $D_n$-modules 
and in that category) quasi-isomorphic to the
complex that is zero except in position $n$ and whose $n$-th entry is
the right $D_n$-module $D_n/(\del_1,\ldots,\del_n)\cdot D_n=\Omega$. A
corresponding statement holds for $\D$.\qed
\end{lem}
\mylabel{idea-O-T}
Notice that 
$\Omega^\bullet(\O_U)$ is the complex
$\Omega^\bullet(\D)\otimes_{\D}\O_U$. 
It follows from the lemma that since
$\Omega^\bullet(D_n)$ is a complex of free $D_n$-modules,
$\Omega^\bullet(\O_U)$ is the complex that
computes ${\mathcal T}or_{n-\bullet}^\D(\Omega(\D),\O_U)$.

The hypercohomology of this complex will simply be the cohomology of
the global section of this complex, because all sheaves in
$\Omega^\bullet(\D)\otimes_{\D}\O_U$ 
are quasi-coherent and $U$ is affine.

So the de Rham cohomology of $U$ is $\tor_{n-\bullet}^{D_n}(\Omega,R_f)$ where 
$Y=\var(f)$.

Let $A^\bullet $ be a resolution of $R_f$ by finitely generated free
$D_n$-modules in
the category of left $D_n$-modules of length greater than $n$. 
That this is possible follows for example
from the fact that $D_n$ is left-Noetherian and that $R_f$
is $D_n$-cyclic (\cite{B}).

Then the cohomology of $\Omega\otimes_{D_n}A^\bullet$ is the de
Rham cohomology of $U$ with coefficients in $\C$ shifted by $n$, since
$H^i(\Omega\otimes_{D_n}A^\bullet)=\tor^{D_n}_{-i}(\Omega,R_f)$ and
$\tor^{D_n}_{-i}(\Omega,R_f)=0$ for $i<0$ and $i>n$.

T.\ Oaku and N.\ Takayama gave an algorithm in \cite{O-T1} for the 
computation of the
cohomology groups of 
this kind of complex. It is in fact explained how one can find the
cohomology groups of the complex 
$D_n/(x_1,\ldots,x_n)\cdot D_n\otimes^L_{D_n} M$
where $M$ is an arbitrary holonomic
$D_n$-module and the tensor product is considered as an element in the
derived category. 
The present problem can be reduced to that case by
applying the Fourier automorphism to $D_n$ which sends $x_i$ to $\del_i$
and $\del_i$ to $-x_i$.

So computation of $H^i_{dR}(U,\C)$ can be summarized as follows
(\cite{O-T1}, algorithm 2.1):
\begin{itemize}
\item Find a suitable finite free resolution $A^\bullet$ of the
$D_n$-module $F(R_f)$, $F(R_f)$ positioned in degree $n$ ($F$ denotes the
Fourier automorphism).  
\item Replace each $D_n$ by the right $D_n$-module 
$\Omega\cong \C[\del_1,\ldots,\del_n]$ in that resolution.
\item Truncate the resolution using the method of \cite{O-T1} 
to a complex of
finite dimensional $\C$-vectorspaces.
\item Take the $i$th cohomology.
\end{itemize}

\subsection{Computing de Rham cohomology for arbitrary $Y$}
\mylabel{MV}
Let $Y$ now be cut out by the $r$ polynomials $f_1,\ldots,f_r$.
The problem arises from the fact that computation of the
hypercohomology of $\Omega^\bullet(\D)\otimes_{\D}\O(U)$ is not just
$\tor_\bullet^{D_n}(\Omega,\Gamma(U,\O(U)))$ anymore, due to the existence 
of higher cohomology of quasi-coherent sheaves on $U$. The strategy is 
to find an open covering of $U$ such that each of the open sets in the
covering is acyclic for cohomology of quasi-coherent sheaves.
\begin{df}
\mylabel{notation-here}
Set $\Re:=$ the ordered nonempty subsets of ${\{1,\ldots,r\}}$. Define
$U_i=X\backslash \var(f_i)$ and more generally for
$I\in\Re$, we define $U_I=\bigcap_{i\in I}U_i$. 

Similarly, set $f_I=\prod_{i\in I}f_i$ with the special cases
$f_I=f_i$ if $I=\{i\}$. Write $\O_{U_I}$ as $\O_I$.
\end{df}
To get started, notice that $U_I=X\backslash \var(f_I)$. This means in
particular, that by Oaku--Takayama the de Rham cohomology groups of
$U_I$ with coefficients in $\C$ are computable as the cohomology of
$\Omega^\bullet(D_n)\otimes_{D_n}R_{f_I}$. 

Notice also that $U=X\backslash Y$ is just the union of all the $U_I$.

In \cite{DRCAV}, page 28, Hartshorne defines how de Rham cohomology of
schemes may be recovered from the de Rham complexes on the open sets
in a finite covering. For our $U$ that works as follows. 

For each $I$ let $X_I=\prod_{i\in I}U_i$. Then $U_I$ embeds in $X_I$
as the diagonal. As $X_I$ is smooth, $\hat\Omega^\bullet(X_I)$
computes de Rham cohomology of $U_I$, the hat denoting completion at
the closed subscheme $U_I\subseteq X_I$.

Consider the direct image $\M_I^\bullet$ of $\hat\Omega^\bullet(X_I)$ in $U$,
induced by the inclusion $j_I:U_I\into U$.

Since $U_I$ is smooth, $\hat\Omega^\bullet(X_I)$  
is
naturally quasi-isomorphic to $\Omega^\bullet(\O_I)$, cf.\
\cite{DRCAV}, Proposition II.1.1.
So $\M_I^\bullet$ is naturally quasi-isomorphic to
$j_{I*}(\Omega^\bullet(\O_I))$, the direct image of
$\Omega^\bullet(\O_I)=\Omega^\bullet(D_n)\otimes_{D_n}\O_I$. 

For $j\not \in I$, the natural maps $X_{I\cup j}\onto X_I$ and
$U_{I\cup j}\into U_I$ give a natural map $
\hat\Omega^\bullet(X_I)\to
\hat\Omega^\bullet(X_{I\cup j})$. Similarly, we get chain maps $\phi_{I,j}: 
\Omega^\bullet(\O_I)\to \Omega^\bullet(\O_{I\cup j})$ induced from the
inclusion $U_{I\cup j}\into U_I$.
It is easy to check that the
natural quasi-isomorphisms from 
$\hat\Omega^\bullet(X_I)$ to $\Omega^\bullet\otimes_{D_n}\O_I$ 
transform the map  $\hat\Omega^\bullet(X_I)\to
\hat\Omega^\bullet(X_{I\cup j})$ into $\phi_{I,j}$. 
So the same is true for the direct images in $U$.

Multiply $\phi_{I,j}$ by $(-1)^{\sgn(I,j)}$, $\sgn(I,j)$
being the number of 
elementary permutations that are needed to make $(I,j)$ an actual element
of $\Re$.

Let us write $\J_I^\bullet:=j_{I*}(\Omega^\bullet(D_n) 
\otimes_{D_n} \O_I)$, a complex of sheaves on $U$.
We will now construct a double complex $\MV(\J)$ out of all the
$\J^\bullet_I$. Let 
$\MV(\J)^{k,l}=\bigoplus_{|I|=l}\J_I^k$. The maps in horizontal ($k$-)
direction
are simply the directs sums of the differentials of the $\J^\bullet_I$
involved, while the vertical ($l$-) maps are defined to be the sums of all
maps which are composed as follows:
$\bigoplus_{|I|=l}\J_I^k\stackrel{nat}{\longrightarrow}
\J_I^k\stackrel{\phi_{I,j}^k}{\longrightarrow}\J_{I\cup j}^k\into 
\bigoplus_{|I'|=l+1}\J_{I'}^k$.  

Notice that this is in fact a double complex (and in particular
anticommutative) due to the sign rule that applies to the
$\phi_{I,j}$.

Then, according to Hartshorne,
 the de Rham cohomology of $U$ is the hypercohomology of the
associated total complex $Tot^\bullet(\MV(\J))$. Of course, 
$\MV(\J)$ is just
the origin of the usual Mayer-Vietoris spectral sequence of de Rham
 cohomology and sometimes called the \v Cech-de Rham complex .  

In the hypersurface case $U$ is affine, so hypercohomology on $U$ was
cohomology of the global sections. Now
we claim
\begin{lem}
The complex $j_{I*}(\Omega^\bullet 
\otimes_{D_n} \O_I)$ consists entirely of sheaves that
have no higher cohomology on $U$.
\proof
In order to see this observe that it is
sufficient to show that $j_{I*}(\O_I)$ has this property, because
$\Omega^i$ is $D_n$-free. If $\E^\bullet_I$ is an injective resolution
of $\O_I$ on $U_I$, then $j_{I*}(\E^\bullet_I)$ is a complex of flasque
sheaves on $U$ as direct images of flasque sheaves are flasque.
Moreover, as $U_I$ is affine, $j_{I*}$ is an exact
functor on quasi-coherent sheaves (and $\O_{U_I}$-morphisms), 
because $R^ij_{I*}(\E^j_I)$ is the sheaf associated to the presheaf 
$V\to H^i(V\cap 
U_I,\E^j_I)$ for open subsets $V$ of $U$. Hence 
we actually get a flasque resolution of
$j_{I*}(\O_I)$. Taking global sections we see that $j_{I*}(\O_I)$ has
no higher cohomology on $U$, as 
$\Gamma(U,j_{I*}(\E^i))=\Gamma(U_I,\E^i)$.\qed 
\end{lem}
\begin{rem}
We note in passing that the proof actually shows that
$H^i(j_{I*}(\F),U)=0$ for positive $i$ and all quasi-coherent $\F$ on
$U_I$.  
\end{rem}
So the complex $Tot^\bullet(\MV(\J))$ consists of
$\Gamma(U,-)$-acyclic sheaves. Thus, in order to compute its
hypercohomology it suffices to compute the
cohomology of the global sections of that complex. We arrive at

\begin{prop}
\mylabel{prop-how-compute-derham}
The de Rham cohomology of of $U$ with coefficients in $\C$, which may be
computed as the hypercohomology of the complex $Tot^\bullet(\MV(\J))$,
agrees with the cohomology of the global sections of
$Tot^\bullet(\MV(\J))$ and can be computed as
$H^\bullet(\Omega^\bullet\otimes_{D_n}MV^\bullet)$, where
\begin{equation}
\label{MV-complex}
MV^\bullet: 0\to
\oplus_{|I|=1}R_{f_I}\to\ldots\to\oplus_{|I|=r}R_{f_I}\to 0.
\end{equation} 
\proof
This follows from the discussion before the proposition, noting that
the global sections on $U$ of  
$j_{I*}(\Omega^\bullet(\O_I))$ are exactly
$\Omega(D_n)\otimes_{D_n}R_{f_I}$ and hence
$\Gamma(U,\MV^\bullet(\J))=\Omega^\bullet(D_n)\otimes_{D_n}MV^\bullet$.
\qed
\end{prop}

For any set of polynomials $\{p_i\}_1^m$, 
the \v Cech complex 
$C^\bullet(R;p_1,\ldots,p_m):=\bigotimes_1^mC^\bullet(p_i)$ is defined by
$C^\bullet(p_i)=
(0\to R\stackrel{1\to\frac{1}{1}}{\longrightarrow} R_{p_i}\to 0)$. 

Notice that $MV^i$ is the $i+1$st entry of the \v Cech 
complex to $f_1,\ldots,f_r$ if $i\geq 0$ and zero otherwise.
\begin{rem}
In the special case where $r=1$, so $I=(f)$, one sees that the complex
$MV^\bullet$ degenerates to $(0\to R_{f_1}\to 0)$ 
reducing to the case from \cite{O-T2}.
\end{rem}
In \cite{W-1}, algorithm 5.1 we gave an algorithm that explicitly
computes the \v Cech complex to a finite set of polynomials as a
complex of finitely generated left $D_n$-modules by means of Gr\"obner
basis computations.

It will now be our task to develop an algorithm that computes the
cohomology of
$\Omega^\bullet\otimes_{D_n}MV^\bullet=\Omega\otimes_{D_n}^LMV^\bullet$.

\section{Computing a certain $D_n$-free complex}
\mylabel{section-free-complex}
In the next two sections we shall be concerned with finding the cohomology of the
complex $\Omega^\bullet\otimes_{D_n} MV^\bullet$ where $MV^\bullet$ is
the Mayer-Vietoris complex (\ref{MV-complex}) to $f_1,\ldots,f_r$, or
more generally the cohomology of $\Omega\otimes^L_{D_n} C^\bullet$
where $C^\bullet$ is an arbitrary complex of finitely generated
$D_n$-modules with specializable cohomology (cf.\ Definition
\ref{V-definitions}).  
In particular, in this section we find a free $D_n$-complex
quasi-isomorphic to a given complex $C^\bullet$ with special properties
related to the so-called $V$-filtration.

We need to introduce some terminology from \cite{O-T1} related to the
$V$-filtration.
\begin{df}
\label{V-definitions}
Fix an integer $d$ with $0\le d\le n$ and set $H=\var(x_1,\ldots,x_d)$. 
For $\alpha\in\Z^n$, we set
$\alpha_H=(\alpha_1,\ldots,\alpha_{d},0,\ldots,0)$. 

On the ring $D_n$ we define the {\em $V_d$-filtration} 
$F_H^j(D_n)$ which consists of all
operators $cx^\alpha\del^\beta$ for which $|\alpha_H|+j\geq|\beta_H|$.
More
generally, on a free $D_n$-module $A=\oplus_1^m D_n\cdot e_i$ we define
$F^j_H(A)[\m]$, where $\m$ is an element of $\Z^m$, as $\sum
F^{j-\m(i)}_H(D_n)\cdot e_i$. We shall call $\m$ the {\em shift vector}.

If $A$ is a free $D_n$-module, the phrase {\em let a shift vector for $A$
be given} will mean the following. First of all we assume that once
and forever a minimal set of generators $\{a_i\}_1^m$ 
for $A$ has been chosen, and
secondly that there is given $\m\in\Z^m$ defining the $V_d$-degree on
$A$ by
the formula of the previous paragraph.

If $M$ is a quotient of the free $D_n$-module $A=\oplus_1^mD_n\cdot
e_i$,  so
$M=A/I$, 
we define the filtration on
$M$ by $F^j_H[\m](M)=F^j_H[\m](A)+I$. 
For submodules $N$ of $A$ we define the $V_d$-filtration by
intersection: $F^j_H[\m](N)=F^j_H[\m](A)\cap N$.

If $A^\bullet$ is a free $D_n$-resolution of the module $M$, $M$ being
positioned in degree zero, we say it
is {\em $V_d$-strict} if there exist shift vectors $\m_i$ such that
$F^j_H[\m_i](A^i)\to F^j_H[\m_{i+1}](A^{i+1})\to F^j_H[\m_{i+2}](A^{i+2})$
is exact for all $i<-1$ and all $j$, and $F^j_H[\m_{-1}](A^{-1})\to
F^j_H[\m_0](A^0)\to F^j_H[\m_{0}](M)\to 0$ is 
exact for all $j$. 

We define the {\em $V_d$-degree} of
an operator, $V_ddeg(P)[\m]$, to be the smallest $i$ such that
$P\in F^i_H[\m](A)$.
\end{df}

It has been shown by T.\ Oaku and N.\ Takayama (\cite{O-T1})
how to compute
for any $D_n$-module $M$ a free $V_d$-strict resolution 
$(A^\bullet[\m_\bullet],\phi^\bullet)$ 
of $M$, $A^i=\oplus_1^{r_i}D_n, r_i=0$ if $i>0$.
The construction given in \cite{O-T1} allows for arbitrary $\m_0$.

The method employed is to construct the free resolution with the usual
technique of finding a Gr\"obner basis for $\ker(A^i\to A^{i+1})$ and
calculating the syzygies on this basis. The trick is  to impose an order that
refines the partial ordering given by $V_d$-degree, together with a
homogenization technique.

The vectors $\m_i$ are obtained for each $A^i$ with falling $i$: 
if $A^i$ maps its generators on a Gr\"obner basis of $\ker(A^{i+1}\to
A^{i+2})$ then the shift component $\m_i(j)$ of the $j$th generator
$e_j$ of
$A^i$ is  defined as $V_ddeg(\phi^i(e_j)[\m_{i+1}])$.

We need to generalize the definitions of \cite{O-T1} to the case where 
the complex $A^\bullet$ is not a resolution.

\begin{df}
A complex of free $D_n$-modules $\ldots\to
A^{k-1}\stackrel{\phi^{k-1}}\to A^k\stackrel{\phi^k}\to 
A^{k+1}\to\cdots$ is said to be {\em $V_d$-adapted at $A^k$ } with respect to
certain shift vectors $\m_{k-1},\m_k,\m_{k+1}$ if
$\phi^kF^j_H[\m_k]A^k\subseteq F^j_H[\m_{k+1}]A^{k+1}$ and also 
$\phi^{k-1}F^j_H[\m_{k-1}]A^{k-1}\subseteq F^j_H[\m_k]A^k$.

We shall say that the complex is {\em $V_d$-strict at $A^k$} if it is
$V_d$-adapted at $A^k$ and moreover $\im\phi^{k-1}\cap
F^j_H[\m_k]A^k=\im(\phi^{k-1}|_{F^j_H[\m_{k-1}]A^{k-1}})$ for all $j$. 

Suppose $P^0$ and $P^1=\oplus_1^mD_n\cdot e_i$
 are free $D_n$-modules, $\phi:P^1\to P^0$ a $D_n$-linear
map and assume that on $P^0$ a shift vector $\m_0$ is given. We define
{\em the obvious shift} on $P^1$ by setting $\m_1(i)=V_ddeg(\phi(e_i))$.

For $1\le d\le n$ we set $\theta_d=x_1\del_1+\ldots+x_d\del_d$. 
Recall that a $D_n$-module $M=A[\m]/I$ 
is called {\em specializable to $H$} if
there is a polynomial in a single variable $b(s)$ such that
$b(\theta_d+j)F^j_H[\m] M\subseteq F^{j-1}_H[\m]M$ for all $j$ (cf.\
\cite{O-T2}). 
Introducing $\gr^j_H[\m]M=(F^j_H[\m]M)\//\/(F^{j-1}_H[\m]M)$, this can be
written as $b(\theta_d+j)\gr^j_H[\m]M=0$.
The polynomial $b$ may depend
on $\m$, while its existence does not (\cite{L-S}).

Notice that independently of $d$, $\gr^\bullet_H(D_n)\cong D_n$, as ring.
\end{df}
The main purpose of this section is to construct for a given finite complex
$0\to C^1\to \ldots\to C^{r-1}\to 0$ with cohomology specializable to $H$ 
a quasi-isomorphic free $V_n$-strict
complex $A^\bullet[\m_\bullet]$. 

\begin{rem}
Notice that if $A^\bullet$ is a free resolution of 
$M$ and $V_d$-strict in our sense
it is also $V_d$-strict in the sense of Oaku/Takayama. In fact our
definition is a natural generalization to complexes that are not
resolutions. 

Moreover, let $\ldots\to
A^{k-1}[\m_{k-1}]\stackrel{\phi^{k-1}}\to A^k[\m_k]\stackrel{\phi^k}\to 
A^{k+1}[\m_{k+1}]\to\cdots$ be $V_d$-strict. 
Then the $V_d$-filtration on $A^k$ induces
a filtration on the $k$-cycles $Z^k=\ker\phi^k$ and since the complex
is $V_d$-strict this gives a natural filtration on the cohomology module
$H^k$, 
\[F^j_H[\m_k]H^k=F^j_H[\m_k]Z^k/\im F^j_H[\m_{k-1}]A^{k-1}.
\] 

Let $B^k$ be $\im\phi^{k-1}$, the $k$-boundaries. 
The short exact sequences $0\to
Z^k\to A^k\to B^{k+1}\to 0$ give rise to short exact sequences of
groups $0\to F^j_H[\m_k]Z^k\to F^j_H[\m_k]A^k\to F^j_H[\m_{k+1}]B^{k+1}\to
0$ since the complex is $V_d$-strict. Similarly, the short exact sequences 
$0\to B^k\to Z^k\to H^k\to 0$ induce short exact sequences $0\to
F^j_H[\m_k]B^k\to F^j_H[\m_k]Z^k\to F^j_H[\m_k]H^k\to 0$.

This in turn induces short exact sequences of the graded objects, 
\[
0\to \gr^j_H[\m_k](B^k)\to \gr^j_H[\m_k](Z^k)\to \gr^j_H[\m_k](H^k)\to 0
\]
and
\[
0\to \gr^j_H[\m_k](Z^k)\to \gr^j_H[\m_k](C^k)\to 
 \gr^j_H[\m_{k+1}](B^{k+1})\to 0.
\]
These sequences are the main feature of $V_d$-strict complexes.
\end{rem}
We shall break the construction of $A^\bullet[\m_\bullet]$ into
several steps.
\begin{lem}
\mylabel{lemma1}
Let $0\to P_A/I_A\to P_B/I_B\to P_C/I_C\to 0$ be exact and assume that 
on $P_C$ there is given a shift vector $\m_C$.

Then $P_B/I_B$ can be replaced by a certain other quotient of a free
module isomorphic to $P_B/I_B$, such that there exist shift vectors
$\m_A,\m_B$ making the sequence $V_d$-strict.
\proof 
We remark 
that making the sequence $V_d$-adapted is trivial (but not good enough).

Set $Q_B=P_A\oplus P_C$ and define $Q_B\to P_B/I_B$ as $Q_B\onto
P_A\onto P_A/I_A\to P_B/I_B$ 
for the $P_A$-part, while for the $P_C$-part we pick some map from
$P_C$ to $P_B/I_B$ that lifts $P_B/I_B\to P_C/I_C$ 
and then map $Q_B\onto P_C\to P_B/I_B$.

Then $Q_B$ projects onto $P_B/I_B$. Let's say the kernel is
$I_{A,C}$. Notice that $I_{A,C}$ contains $I_A\oplus 0$, which
corresponds to
the injection $P_A/I_A\into Q_B/I_{A,C}$. 

Now define the shift on $Q_B$ by taking the given shift from $P_C$ on the
second component, and for the generators of $P_A$ take any shift of
your choice.

It is clear that the resulting short exact sequence $0\to P_A/I_A\to
Q_B/I_{A,C} \to P_C\to 0$ is $V_d$-adapted. It is just as clear that it is 
strict at $P_C/I_C$ and $P_A/I_A$. Let $b=\sum
\alpha_ia_i+\sum\gamma_ic_i$ be an element of $Q_B$ that is
sent to zero in $P_C/I_C$. That means that $\sum \gamma_ic_i\in
I_C$. Since $I_{A,C}$ contains for all elements $c\in I_C$ an element
$(a_c,c)$ (after all, modulo the first component of $Q_B$,
$0\oplus I_C$-elements must be zero!) we can transform our element $b$
into an element of $I_A\oplus 0$. 
We would like this element to have $V_d$-degree at most 
$V_ddeg(b)$. 

Remember that we can still choose the shifts for $P_A$! Take a
$V_d$-strict G-basis for $I_C$ with respect to the given shift $\m_C$. Let it
have elements $\{n_i\}$. For each $n_i$ find a $m_i$ in $P_A$ such that
$(m_i,n_i)\in I_{A,C}$. Now define the shift in $P_A$ in such a way
that $m_i$ has $V_d$-degree at most equal to the $V_d$-degree of $n_i$
for all $i$. 
We aim to show that with this shift the sequence is $V_d$-strict.

Return to $b=\sum \alpha_ia_i+\sum\gamma_ic_i$. The $V_d$-degree of $b$
is the maximum of the degrees of the two sums. Since $\sum\gamma_ic_i$ 
is in $I_C$, we can write it as a sum $\sum\delta_in_i$, where the
$V_d$-degree of the sum, let's call it $e$, is the $V_d$-degree of the largest
summand in the sum, because the $n_i$ form a G-basis.

Modulo $I_{A,C}$, this is the same as the sum $-\sum\delta_im_i$,
which has lower or equal $V_d$-degree, by construction of the shift on
$P_A$. 

Then $\sum \alpha_ia_i-\sum\delta_im_i\in P_A$ is an expression that maps
onto $b$, modulo $I_{A,C}$, 
and has degree at most equal to $e$. Strictness follows.

We are done.
\qed
\end{lem}
Notice that this creates a commutative diagram with exact and
$V_d$-strict rows and
columns 
\[
\diagram
&0&0&0&\\
0\rto&P_A/I_A\uto\rto^{\phi_A}&P_B/I_B\uto\rto^{\phi_B}&P_C/I_C\uto\rto&0\\
0\rto&P_A\uto\rto&Q_B\uto\rto&P_C\uto\rto&0\\
0\rto&I_A\uto\rto&I_{A,C}\uto\rto&I_C\uto\rto&0\\
&0\uto&0\uto&0\uto&
\enddiagram
\]
We will have need of a slight improvement of lemma \ref{lemma1}:
\begin{lem}
\mylabel{lemma2}
Suppose we have 2 short exact sequences 
\[
0\to P_A/I_A\to P_B/I_B\to P_C/I_C\to 0
\]
and
\[
0\to P_{D}/I_{D}\to P_A/I_A\to P_F/I_F\to 0
\]
and assume that on $P_C$ is given a shift vector $\m_C$. Then one can
rewrite 
$P_A/I_A$ as $Q_A/I_{D,F}$
and 
$P_B/I_B$ as $Q_B/I_{D,F,C}$ 
and find shift vectors $\m_A, \m_B, 
\m_{D}, \m_F$ such that the resulting 2 sequences are exact and
$V_d$-strict.
\proof
First use the first half of the proof of the previous lemma to write
$P_A/I_A$ as $Q_A/I_{D,F}$ and then with that representation of
$P_A/I_A$ rewrite $P_B/I_B$ as $Q_B/I_{D,F,C}$. So in particular,
$Q_B=Q_A\oplus P_C=P_{D}\oplus P_F\oplus P_C$.

In order to find the proper shift vectors, proceed as follows:
\begin{enumerate}
\item Take a G-basis $\{c_i\}$ for $I_C$ with respect to an order
refining the 
$V_d$-filtration on $P_C$ relative to the given shift $\m_C$. For all $i$
find $a_i=(d_i',f_i')$ such that $(d_i',f_i',c_i)\in I_{D,F,C}$.
\item Pick a shift on $F$ such that $V_ddeg(f_i')\le V_ddeg(c_i)$ for all
$i$.
\item Compute a G-basis $\{f_i\}$ of $I_F$ using an order that refines
$V_d$-degree on $F$, using the shift we just found. For all $i$ find $d_i$ 
with $(d_i,f_i)\in I_{D,F}$.
\item Pick a shift on $P_{D}$ such that $V_ddeg(d_i)\le V_ddeg(f_i)$ for all 
$i$ and $V_ddeg(d_i')\le V_ddeg(c_i)$ for all $i$.
\end{enumerate}
By arguments similar to those that prove lemma \ref{lemma1}, the
sequences are $V_d$-strict.
\qed
\end{lem}
Lemma \ref{lemma1} and \ref{lemma2} providing the basis for the
construction, the following result is the inductive step:
\begin{lem}
\mylabel{lemma3}
Let $A,B,C$ be three submodules of free modules $F_A,F_B,F_C$. Assume
that $0\to A\to B\to C\to 0$ is exact and $V_d$-strict, relative to some
shift vectors on $F_A,F_B,F_C$.

Then one can construct a diagram
\begin{eqnarray}
\diagram
&0&0&0&\\
0\rto&A\uto\rto^{\phi_A}&B\uto\rto^{\phi_B}&C\uto\rto&0\\
0\rto&P_A\uto\rto&P_B\uto\rto&P_C\uto\rto&0\\
0\rto&K_A\uto\rto&K_B\uto\rto&K_C\uto\rto&0\\
&0\uto&0\uto&0\uto&
\enddiagram
\label{diagram-type}
\end{eqnarray}

such that 
\begin{itemize}
\item all $P_X$ are free,
\item all rows and columns are exact,
\item there are shift vectors $\m_A,\m_B,\m_C$ such that if
$P_A,P_B,P_C$ are shifted accordingly, all rows and columns become
$V_d$-strict.
\end{itemize}
\proof
Let $\{a_i\}, \{b_i\}$ and $\{c_i\}\cup\phi_B(\{b_i\})$  
be G-bases for $A,B,C$ with respect to 
an order on $F_A,F_B,F_C$ that refines $V_d$-degree.

Let $P_A$ be a free module on the symbols $\{e_{a_i}\}$, define a 
projection from $P_A$ to $A$ in the obvious way. 
Let $P_C$ be a free module on the symbols $\{e_{b_i}\}\cup
\{e_{c_i}\}$. 
Define the projection $P_C\to C$ by $e_{c_i}\to c_i, e_{b_i}\to
\phi_B(b_i)$.
Define degree shifts on $P_A,P_C$ in the obvious way.

Set $P_B=P_A\oplus P_C$. Define $P_B\to B$ on the $P_A$-part in the
obvious way, and similarly on the $P_C$-part that corresponds to
$\{b_i\}$. For the generators corresponding to $\{c_i\}$, use a lift
$\psi: P_C\to B$ for $\phi_B$ 
(which exists as $P_C$ is free) that satisfies:
$\psi(e_{c_i})$ is of no higher (shifted) $V_d$-degree than $V_ddeg(c_i)$
(which exists because $0\to A\to B\to C\to 0$ is $V_d$-strict).
Set $K_A,K_B,K_C$ to be the corresponding kernels.

It is clear that $0\to P_A\to P_B$ and $P_B\to
P_C\to 0$ are $V_d$-strict. If an element of $V_d$-degree $e$ is in
the kernel of $P_B\to P_C$, then its second  component (the one in
$P_C$) is zero, 
so the $V_d$-degree came from the $P_A$-component. Hence the second row is 
$V_d$-strict.
Then automatically the third row is too.

By \cite{O-T2}, the remarks after proposition 3.11, 
the outer columns are $V_d$-strict.

Let $b\in B$ be of $V_d$-degree $d$. Then $b=\sum \alpha_ib_i$ where
$V_ddeg(\alpha_i)+V_ddeg(b_i)\le V_ddeg(b)$ for all $i$, since $\{b_i\}$ is
a G-basis and the image of the element $\sum \alpha_ie_{b_i}\in P_B$
is $b$. Moreover, by our definitions, the $V_d$-degree of this sum 
in $P_B$ is at most the degree of the image in $B$, which is $e$. It
follows that $P_B\to B\to 0$ is $V_d$-strict, and hence the whole column.
\qed
\end{lem}

We need a version of lemma \ref{lemma3} of the type of \ref{lemma2}.

\begin{lem}
\mylabel{lemma4}
Assume we have 2 exact $V_d$-strict sequences $0\to A\to B\to C\to 0$ and
$0\to D\to A\stackrel{\phi}{\to} F\to 0$ of submodules of free
modules. 
Suppose moreover
that for $D$ we are given a 
short exact $V_d$-strict sequence $0\to K_{D}\to 
P_{D}[\m_{D}]\to D\to 0$ where the free module 
$P_{D}$ maps its generators on a Gr\"obner
basis of $D$ with respect to the given shift vector. Then for
$X\in\{A,B,C,F\}$ one can find free module $P_X$ and  
submodules 
$K_X$ creating commutative diagrams with $V_d$-strict and exact rows and
columns of the type (\ref{diagram-type}) over both given short exact
sequences. 
\proof
Let $\{a_i\}$ be a G-basis for $A$ and assume $\{\phi(a_i)\}\cup
\{f_i\}$ 
is one for $F$. Then set
$P_F$ to be the free $D_n$-module on the symbols $\{\phi(a_i)\}\cup \{f_i\}$
and $P_A=P_{D}\oplus P_F$. Define maps $P_F\to F$ in the obvious way and 
$P_A\to A$ as $(P_A\onto P_{D}\to D\to A) + (P_A\to P_F\to A)$ where the 
map from $P_F\to A$ is defined by $e_{\phi(a_i)}\to a_i$ and
$e_{f_i}\to g_i$ where $g_i$ is a preimage of $f_i$ in $A$ of lesser
or equal $V_d$-degree as $f_i$.

Let $K_A, K_F$ be the kernels.

This gives as in lemma \ref{lemma3} a $V_d$-strict exact commutative
diagram over $0\to A\to D\to F\to 0$.

Now use the sequence $0\to K_A\to P_A\to A\to 0$ to build a diagram
over $0\to A\to B\to C\to 0$.\qed
\end{lem}

We have now assembled enough machinery to to find for all complexes of 
$D_n$-modules the cohomology of which is 
specializable to $H=\var(x_1,\ldots,x_d)$ 
a quasi-isomorphic $V_d$-strict
complex. 

So suppose we are given such  a complex 
$0\to C^1\to \ldots\to C^{r-1}\to 0$. As
a first step, write down all the short exact sequences $0\to B^k\to
Z^k\to H^k\to 0$ and $0\to Z^k\to C^k\to B^{k+1}\to 0$. That is to
say, find representations of these modules and maps in terms of finitely
generated free modules modulo a finite number of relations.
Observe that all $B^{i+1}, C^{i+1}, Z^{i+1}, H^{i+1}$ are zero for $i\geq r-1$.

Invoke lemma \ref{lemma2} to find a
presentation for $Z^{r-1}$ and for $C^{r-1}$ 
together with shift vectors $\m_{Z,r},
\m_{C,r}, \m_{H,r}$ 
and $\m_{B,r}$ such that there are commutative diagrams
\[
\diagram
&0&0&0&\\
0\rto&Z^{r-1}\uto\rto&C^{r-1}\uto\rto&B^{r}=0\uto\rto&0\\
0\rto&P_{Z,r-1}[\m_{Z,r-1}]\uto\rto
 &P_{C,r-1}[\m_{C,r-1}]\uto\rto&P_{B,r}[\m_{B,r}]=0\uto\rto&0\\ 
0\rto&K_{Z,r-1}[\m_{Z,r-1}]\uto\rto
&K_{C,r-1}[\m_{C,r-1}]\uto\rto&K_{B,r}[\m_{B,r}]=0\uto\rto&0\\ 
&0\uto&0\uto&0\uto&
\enddiagram
\]
and 
\[
\diagram
&0&0&0&\\
0\rto&B^{r-1}\uto\rto&Z^{r-1}\uto\rto&H^{r-1}\uto\rto&0\\
0\rto&P_{B,r-1}[\m_{B,r-1}]\uto\rto
 &P_{Z,r-1}[\m_{Z,r-1}]\uto\rto&P_{H,r-1}[\m_{H,r-1}]\uto\rto&0\\ 
0\rto&K_{B,r-1}[\m_{B,r-1}]\uto\rto
&K_{Z,r-1}[\m_{Z,r-1}]\uto\rto&K_{H,r-1}[\m_{H,r-1}]\uto\rto&0\\ 
&0\uto&0\uto&0\uto&
\enddiagram
\]

with all $P_X$ free, that has exact and $V_d$-strict rows and columns.

Then invoke the lemma again, this time starting with the shift just
obtained on
$B^{r-1}$ and constructing representations for $Z^{r-2}, C^{r-2},
$ and shifts on $P_{Z,r-2}, P_{C,r-2}, P_{H,r-2}, P_{B,r-2}$.

And so on. Repetition leads to $V_d$-strict commutative diagrams
with exact rows and columns 
\[
\diagram
&0&0&0&\\
0\rto&Z^i\uto\rto&C^i\uto\rto&B^{i+1}\uto\rto&0\\
0\rto&P_{Z,i}[\m_{Z,i}]\uto\rto
 &P_{C,i}[\m_{C,i}]\uto\rto&P_{B,i+1}[\m_{B,i+1}]\uto\rto&0\\ 
0\rto&K_{Z,i}[\m_{Z,i}]\uto\rto
&K_{C,i}[\m_{C,i}]\uto\rto&K_{B,i+1}[\m_{B,i+1}]\uto\rto&0\\ 
&0\uto&0\uto&0\uto&
\enddiagram
\]
and
\[
\diagram
&0&0&0&\\
0\rto&B^i\uto\rto&Z^i\uto\rto&H^{i}\uto\rto&0\\
0\rto&P_{B,i}[\m_{B,i}]\uto\rto
 &P_{Z,i}[\m_{Z,i}]\uto\rto&P_{H,i}[\m_{H,i}]\uto\rto&0\\ 
0\rto&K_{B,i}[\m_{B,i}]\uto\rto
&K_{Z,i}[\m_{Z,i}]\uto\rto&K_{H,i}[\m_{H,i}]\uto\rto&0\\ 
&0\uto&0\uto&0\uto&
\enddiagram
\]
for $0\le i< r$. The point of this procedure is the creation of a
presentation of $C^i$ as $P_{C,i}/K_{C,i}$ with $V_d$-strict maps
between these modules.

Now we assemble a resolution for $C^\bullet$ as follows. 
First find a $V_d$-strict
resolution using the method of \cite{O-T2} for  $B^0$. 
With lemma \ref{lemma4} find a resolution for $Z^0, H^0, C^0, B^1$.
Then feed the obtained resolution for $B^1$ into lemma \ref{lemma4},
resulting in resolutions for $Z^1, H^1, C^1, B^2$. Et cetera.

Let us denote the $k$-th module of the resolution for $X^i$ (with $X$
being $Z,C,B$ or $H$) by $P^k_{X,i}$. We define a map $\delta^k_{C,i}$
from $P^k_{C,i}$
to $P^k_{C,i+1}$ as the combined maps $P^k_{C,i}\to P^k_{B,i+1} 
\to P^k_{Z,i+1} \to P^k_{C,i+1}$, 
multiplyied by $(-1)^{k}$. So up to sign $\delta_{C,i}^k$ 
 is $P^k_{C,i}=P^k_{B,i}\oplus
P^k_{H,i}\oplus P^k_{B,i-1}\onto P^k_{B,i+1}\into P^k_{B,i+1}\oplus
P^k_{H,i+1}\oplus P^k_{B,i+1}$. 
Clearly $P^k_{C,i}\to P^k_{C,i+1}\to
P^k_{C,i+2}$ is the zero map and the horizontal maps (in $i$-direction)
are $V_d$-strict.

We have created a double complex $P^\bullet_{C,\bullet}$ of free
$D_n$-modules. Moreover, the associated total complex is
quasi-isomorphic to $C^\bullet$ and clearly $V_d$-adapted.

\begin{prop}
$Tot^\bullet(P^\bullet_{C,\bullet})$ is in fact $V_d$-strict. 
\proof
To that end assume
that the element $p=p^0_i\oplus p^1_{i+1}\oplus\cdots\oplus
p^{r-1-i}_{r-1}\in Tot^i(P^\bullet_{C,\bullet})= P^0_{C,i}\oplus
P^1_{C,i+1}\oplus\cdots \oplus P^{r-1-i}_{C,r-1}$
is in the image of 
the total differential $\delta_{T}$, and that the $V_d$-degree of $p$ under
the shift vectors is $e$. 
We need to take a closer look at the maps and
modules in front of us.

$P^k_{C,i}$ is by construction $P^k_{B,i}\oplus P^k_{H,i}\oplus
P^k_{B,i+1}$. 
The map
$P^{k+1}_{C,i}=P^{k+1}_{B,i}\oplus P^{k+1}_{H,i}\oplus 
P^{k+1}_{B,i+1}\to P^k_{C,i}=P^k_{B,i}\oplus P^k_{H,i}\oplus
P^k_{B,i+1}$ is on the first component the differential from the
resolution $P^\bullet_{B,i}$ while the map from $P^{k+1}_{H,i}\oplus 
P^{k+1}_{B,i+1}\to P^k_{C,i}$ is defined using certain
lifts, obtained while using lemma \ref{lemma3}. 
Inspection shows that the matrix which represents $P^{k+1}_{C,i}\to
P^k_{C,i}$  looks 
like this:
\[
\left(
\begin{array}{ccc}
\delta_{B,i}^{k+1}&\phi_1^{k+1}&\psi_1^{k+1}\\
0&\delta_{H,i}^{k+1}&\psi_2^{k+1}\\
0&0&\delta_{B,i+1}^{k+1}
\end{array}
\right)
\]
where the $\delta$ are the differentials of the various resolutions
for $B$ and $H$ and $\phi_1^{k+1}, \psi_1^{k+1},\psi_2^{k+1}$ 
are the maps that are used 
to produce the mentioned lifts. Note that $\psi_1,\psi_2,\phi$ are all
$V_d$-adapted by construction. 

We shall argue by falling induction on the variable $s$, starting with
$s=r-1$, that the components $p^{s-i}_{s}$ of $p$ may
be assumed to 
be zero modulo images of degree no greater than $e$ under the total
differential. 
We will at the same time show that we may
assume that the third component of $p^{s-i}_s$ is zero. For $s<r-1$
this will follow from the induction. For $s=r-1$ it follows from the
fact that $B^{r}=0$.

So assume that $0\le s\le r-1$, that $p$ has only zero components
beyond the $s$-th component and that the third piece (to $P_{B,s+1}$)
of the $s$-th component of $p$ is zero.

The following lemma will essentially show that our $p$ is then 
in fact image of an
element in $Tot^{i+1}(P^\bullet_{C,\bullet})$ 
with zero component in $P^{s-i+1}_{B,s+1}$ and
only zeros in all columns beyond the $s$th. 

\begin{lem}
Let $(a,b,0)\in P^{s-i}_{C,s}$ and assume
$(a,b,0)=\delta_C(\alpha,\beta,\gamma)$ with $(\alpha,\beta,\gamma)\in
P^{s-i+1}_{C,s}$. Then $(\alpha,\beta,0)=\delta_C(\alpha',\beta',0)$ for
some $(\alpha',\beta',0)\in P^{s-i+1}_{C,s}$
where
$V_ddeg(\alpha',\beta')\le V_ddeg(\alpha,\beta,\gamma)$.
\proof 
By construction 
$(\psi_1(\gamma),\psi_2(\gamma))$ is in $\ker(P^{s-i}_{B,s}\oplus
P^{s-i}_{H,s}\to P^{s-i-1}_{B,s}\oplus
P^{s-i-1}_{H,s})$. Since this kernel is exactly
$\delta_Z(P^{s-i+1}_{B,s}\oplus P^{s-i+1}_{H,s})$,
$(\psi_1(\gamma),\psi_2(\gamma))=\delta_Z(\alpha'',\beta'')$ where we
can pick $\alpha''$ and $\beta''$ to be of $V_d$-degree at most
$V_ddeg(\gamma)$. Thus
$\delta_C(\alpha,\beta,\gamma)=
(\delta_B\alpha+\psi_1\gamma+\phi\beta,\delta_H\beta+\psi_2\gamma,0)=
\delta_C(\alpha,\beta,0)+\delta_C(\alpha'',\beta'',0)$.
Set $(\alpha',\beta')=(\alpha''+\alpha,\beta''+\beta)$.\qed
\end{lem}
Now write $p^{s-i}_s=(a,b,0)$. The lemma tells us that $b$
equals 
$\delta_H(b_1)$ for some $b_1\in P^{s-i+1}_{H,s}$, of $V_d$-degree at most
$e$ because $P^\bullet_{H,s}$ is $V_d$-strict. What can we say about 
$p-\delta_T(b_1)$, which we call $p$ from now on?

Certainly the $V_d$-degree is at most $e$, it is in the image of
$\delta_T$, only the first component of $p^{s-i}_s$ is nonzero and all
components beyond the $s$th one are zero. 
\begin{lem}
\mylabel{lemma6}
Suppose $(a,0,0)$ is the image under $\delta_C$ of $(\alpha,\beta,\gamma)\in
P^{s-i+1}_{C,s}$. Then there is $\alpha'\in P^{s-i+1}_{B,s}$ with
$\delta_C(\alpha',0,0)=(a,0,0)$ and $\alpha'$ can be chosen to be of
$V_d$-degree no bigger than $V_ddeg(a)$. 
\proof
By the previous lemma, we can assume that $(a,0,0)$ is the image of
$(\alpha,\beta,0)$. By construction, $\phi(b)\in\ker(P^{s-i}_{B,s}\to
P^{s-i-1}_{B,s})=\im(P^{s-i+1}_{B,s}\to P^{s-i}_{B,s})$. 
As $P^\bullet_{B,s}$ is
a $V_d$-strict resolution, $\phi(b)=\delta_B(\alpha'')$ for some $\alpha''\in
P^{s-i+1}_{B,s}$. Hence
$\delta_C(\alpha,\beta,0)=
(\delta_B\alpha+\phi\beta,0,0)=\delta_C(\alpha,0,0)+\delta_C(\alpha'',0,0)$.

Since $\delta_B$ is $V_d$-strict and $\phi$ is $V_d$-adapted, 
we can choose $\alpha'=\alpha+\alpha''$
to be of $V_d$-degree at most $V_ddeg(a)$.\qed
\end{lem}
The component of $p$ in $P^{s-i}_{C,s}$ looks like $(a,0,0)$. Lemma
 \ref{lemma6} tells us that since $p$ is
an image under $\delta_T$, 
$a=\delta_B(\alpha)+(-1)^{s-1}\alpha'$ where $\alpha$ lives in
$P^{s-i+1}_{B,s}$ (a component of $P^{s-i+1}_{C,s}$) 
and $\alpha'\in P^{s-i}_{B,s}$ (a component of $P^{s-i}_{C,s-1}$).

Then the third component of $p$ in $P^{s-i-1}_{C,s-1}$ must be exactly
$\delta_B(\alpha')=(-1)^{s-1}\delta_B(a)$. 
Replace $p$ by $p-\delta_T((-1)^sa)$, this $-a$
positioned in the $P^{s-i-1}_{B,s}$-component of $P^{s-i-1}_{C,s-1}$.

The result is a $p$ with zero component for $P^{s-i}_{C,s}$ and zero
component in the third component for $P^{s-i-1}_{C,s-1}$ which differs
from the original $p$ by the $\delta_T$-boundary of an element of
$V_d$-degree at most $e$.

Now repeat the argument for $p^{r-i-2}_{r-2}, p^{r-i-3}_{r-3},\ldots$.\qed
\end{prop}

We have proved 
\begin{thm}
If $C^\bullet$ is a complex of $D_n$-modules that is bounded below and
above and presentations of all $C^i$ in terms of generators and
relations are given, then one can produce a $V_d$-strict complex of free
$D_n$-modules that is quasi-isomorphic to $C^\bullet$.
\qed
\end{thm}
\begin{rem}
It is not true that the total complex associated to any double complex
with $V_d$-strict rows and columns is $V_d$-strict. This would be equivalent
to saying that all finite subsets of a free $D_n$-module form a G-basis
for any order refining $V_d$-degree. Consider for example the diagram
$\diagram 
D_1[1]\rto^{\cdot\del_1}&D_1[0]\\
0\rto\uto&D_1[1]\uto^{\cdot(\del_1-1)}\enddiagram$. Here, $1=1\cdot
\del_1-1\cdot(\del_1-1)
\in
\im(F^1_H[1,1](Tot^1))$ but it is not in $\im(F^0_H[1,1](Tot^1)\to
F^0_H([0](Tot^0)))$, 
although it is of $V_1$-degree 0.
\end{rem}

\section{Computing Cohomology of $\Omega^\bullet(D_n)\otimes_{D_n} MV^\bullet$}
\mylabel{restrict-complex}
For this section let $A^\bullet[\m_\bullet]$ be 
a
given $V_d$-strict free complex. We shall assume further that the cohomology
modules of $A^\bullet[\m_\bullet]$ are specializable to
$H=\var(x_1,\ldots,x_d)$. 
The main purpose of this section is to determine a suitable truncation
of $A^\bullet[\m_\bullet]$ 
such that the cohomology of
$\Omega\otimes_{D_n} 
A^\bullet$ is captured by the ``tensor product'' of
$\Omega$ with that
truncation. 

Recall that in \cite{O-T1} $A^\bullet$ is a $V_n$-strict resolution of a
specializable module $M$ and the truncation is determined by
considering roots of the $b$-function $b(M)$ corresponding to
restriction to the origin (\cite{O-T1},
Algorithm 5.4). 

\subsection{}
\mylabel{subsection-global-b-function}
Let $H$ be the subspace defined by $x_1=\ldots=x_d=0$, and $A^\bullet$
a $V_d$-strict complex of free $D_n$-modules. 
As a first step,  pick generators $\kappa_{i,l}$ for the kernel modules
$Z^i=Z^i(A^\bullet)$ of
$A^\bullet$. To each of them is associated a degree in the shifted
$V_d$-filtration from $A^i$ which we shall call $\lambda_{i,l}$. 

Let a bar denote cosets of elements of $Z^i$ in $H^i=H^i(A^\bullet)$.
Recall that we agreed to write $\theta_j=x_1\del_1+\ldots+x_j\del_j$
for $1\le j\le n$.  
Since $D_n\cdot \bar\kappa_{i,l}$ is a specializable $D_n$-module, 
there is a $b$-function
$b_{i,l}(\theta_d)$ associated to it which corresponds to the
restriction of 
$D_n\cdot \bar\kappa_{i,l}$ to $x_1=\ldots=x_d=0$. Therefore,
$b_{i,l}(\theta_d) \kappa_{i,l}\in  
F_H^{-1}(D_n)\cdot \kappa_{i,l}+\im (A^{i-1}\to A^i)$. 

Let $b(\theta_d)$ 
be the least common multiple of all $b_{i,l}(\theta_d-\lambda_{i,l})$.

\subsection{}
Then consider the
associated complex of graded $\gr^\bullet_H(D_n)$-modules $\oplus F_H^j
A^\bullet/F_H^{j-1}A^\bullet$. 

Now assume that the $\kappa_{i,l}$ form a G-basis for $Z^i$  under 
the order on $A^i$. 
Then for all $\zeta\in Z^i$, $\zeta=\sum\alpha_{i,l}(\zeta)\kappa_{i,l}$
with $V_ddeg(\alpha_{i,l}(\zeta)\kappa_{i,l})\le V_ddeg(\zeta)$. Hence 
the $\kappa_{i,l}$ are generators for $\gr^\bullet_H(H^i)$  
and moreover 
$\gr_H^j(Z^i)=\sum \gr_H^{j-\lambda_{i,l}}(D_n)\kappa_{i,l}$. Since
$\im(F_H^jA^{i-1}[\m_{i-1}]\to F^j_HA^i[\m_i]) =F^j_HA^i[\m_i]\cap
\im(A_H^{i-1}\to A^i)$ we have  $\gr_H^j(H^iA^\bullet)=\sum
\gr^{j-\lambda_{i,l}}(D_n)\bar\kappa_{i,l} $.

Then observe the following:
\begin{eqnarray}
\label{b-function-kills}
b(\theta_d+j)\gr_H^j(H^iA^\bullet)&=&b(\theta_d+j)\sum
\gr_H^{j-\lambda_{i,l}}(D_n)\bar\kappa_{i,l} \\
&=&\sum \gr_H^{j-\lambda_{i,l}} 
(D_n)b(\theta_d+\lambda_{i,l})\bar\kappa_{i,l}\\
&=&0
\end{eqnarray}
because $b(\theta_d+\lambda_{i,l})$ sends $\kappa_{i,l}$ into
$F_H^{-1}(D_n)\cdot \kappa_{i,l}+\im(A^{i-1}\to A^i)$, which is zero in
$\gr_H^{\lambda_{i,l}}[\m_i](H^i(A^\bullet))$.  

\begin{rem}
$b(\theta_d)$ is then a multiple of the $b$-function
of $H^i(A^\bullet)$ with respect to the given shift vectors.
\end{rem}

This paves the way for a result related to \cite{O-T1}, Proposition
5.2. The proof is very similar to the one given there. 

We need to introduce a number of Koszul complexes.
Let $\L$ be a $\gr_H^\bullet(D_n)$-module and let $\L_j, j\in\Z$ be subgroups
of $\L$ such that $x_i\L_j\subseteq \L_{j-1}$ for $1\le i\le d$. 
In that case we will
say that the $\L_i$ give an  {\em $H$-filtration} for $\L$. 
For any integer $k$ let 
$\K^\bullet (\L,x_1,\ldots,x_d)[k]$ be the Koszul complex 
\[
0\to\L_{k+d}\otimes_\Z \bigwedge^0\Z^{d}\to\L_{k+d-1}\otimes_\Z
\bigwedge^1\Z^{d}\to
\cdots\to \L_k\otimes_\Z \bigwedge^{d}\Z^{d}\to 0
\]
equipped with the usual Koszul maps $\delta(u\otimes
e_{i_1}\wedge\cdots\wedge e_{i_j})=\sum_l x_lu\otimes
e_l\wedge e_{i_1}\wedge\cdots\wedge e_{i_j}$. 

Unifying all the graded pieces, let $\K^\bullet(\L_\bullet,x_1,\ldots,x_d)$ 
be the usual Koszul complex of $\L$ relative to $x_1,\ldots,x_d$.

More generally, for a complex of $H$-filtered $\gr^\bullet_H(D_n)$-modules
$(\L^\bullet,\delta^\bullet)$ 
with morphisms that respect the $H$-filtration we define 
inductively $\K^\bullet(\L^\bullet[k],x_1,\ldots,x_d)$ as the total
complex of the double complex
\[
\diagram
\K^\bullet (\L^{\bullet}_{\bullet},x_1,\ldots,x_{d-1})[k]\\
\K^\bullet (\L^{\bullet}_{\bullet+1},x_1,\ldots,x_{d-1})\uto^{(-1)^ix_d}[k]\\
\enddiagram
\]
where
$\K^\bullet(\L^\bullet_\bullet,\emptyset)[k]=
\L^\bullet_\bullet[k]=(\cdots\to
\L^i_k\to \L^{i+1}_k\to\cdots)$, the $k$-th piece of the original
complex. 

Notice that $\K^\bullet(\L^\bullet,x_1,\ldots,x_d)[k]$ is the graded
component of the usual Koszul complex 
$\K^\bullet(\L^\bullet,x_1,\ldots,x_d)$ 
associated to $\L^\bullet $ and $x_1,\ldots,x_d$ 
that ``ends'' in $V_d$-degree $k$.

The following theorem explains which graded pieces
$\K^\bullet(\L_\bullet^\bullet,x_1,\ldots,x_d)[k]$ of
$\K^\bullet(\L_\bullet^\bullet,x_1,\ldots,x_d)$ are responsible for
nontrivial cohomology pieces of 
$\K^\bullet(\L_\bullet^\bullet,x_1,\ldots,x_d)$.

\begin{thm}
\mylabel{b-function-complex}
Suppose given is a complex of graded $\gr_H^\bullet(D_n)$-modules
$\L^\bullet$ (i.e., $\gr_H^k(D_n)\L^\bullet_j\subseteq \L^\bullet_{j+k}$)
where the maps between the $\L^i$ preserve the grading. Assume that
there is a polynomial $b(\theta)$ in $\C[\theta]$ that satisfies
$b(\theta_d+j)\ker(\L^i_j\to \L^{i+1}_j)\subseteq \im(\L^{i-1}_j\to
\L^i_j)$
for all $j$ and all $i$. Let $k$ be an integer for which
$b(k)\not=0$. Then $\K^\bullet(\L^\bullet,x_1,\ldots,x_d)[k]$ is exact.
\proof
The complex
$\K^\bullet(\L^\bullet,x_d)$ is $V_{d-1}$-graded in the obvious way.  
Here is the essential idea of the argument:
\par{Induction claim. } If $b(\theta_d+j)$
kills cohomology in $\L^\bullet$ of $V_d$-degree $j$, 
then $b^2(\theta_{d-1}+j)$
kills cohomology of $V_{d-1}$-degree $j$ in $\K^\bullet(\L^\bullet,x_d)$. In
other words, cohomology of $\K^\bullet(\L^\bullet,x_d)[j]$.

We may assume that $b(\theta_d)$ is not a constant since otherwise
$\L^\bullet$ is exact and a spectral sequence argument shows that then
$\K^\bullet(\L^\bullet,x_d)$ is exact as well.

So assume that $(u^{i+1}_{j+1},u^i_j)$ is in
$\K^\bullet(\L^\bullet,x_d)[j]$
(whose $i$-th  piece is 
$\L^{i+1}_j\oplus \L^i_j$) and suppose this element is in the kernel of
the differential in  $\K^\bullet(\L^\bullet,x_d)$. Then we must have
\begin{eqnarray}
\delta^{i+1}u^{i+1}_{j+1}=0, \label{first}\\
x_du^{i+1}_{j+1}+\delta^iu^i_j=0, \label{second}
\end{eqnarray}
$\delta$ denoting the boundary map in $\L^\bullet$.
By hypothesis on $b$, $b(\theta_d+j+1)u^{i+1}_{j+1}=\delta^iu^i_{j+1}$ 
for some $u^i_{j+1}\in \L^i_{j+1}$. So
$b(\theta_{d-1}+j+\del_dx_d)u^{i+1}_{j+1}=\delta^iu^i_{j+1}$ and therefore 
$b(\theta_{d-1}+j)u^{i+1}_{j+1}+\del_dPx_du^{i+1}_{j+1}=\delta^iu^i_{j+1}$ 
for some $V_d$-homogeneous $P\in F^0_H(D_n)\backslash F^{-1}_H(D_n)$.
Hence
$b(\theta_{d-1}+j)u^{i+1}_{j+1}=\delta^i(u^i_{j+1}-\del_dPu^i_j)$ 
using relation (\ref{second}). Let us write this as
\begin{eqnarray}
b(\theta_{d-1}+j)u^{i+1}_{j+1}=\delta^i(a^i_{j+1}),\label{third}
\end{eqnarray}
$a^i_{j+1}\in(\L^i_{j+1}\oplus 0)\subset
\K^{i-1}(\L^\bullet,x_d)[j+1]$. 

This implies that if $(u^{i+1}_{j+1},u^i_j)$ is in the kernel of
$\delta_T$, the differential on $\K^\bullet(\L^\bullet,x_d)$, then
$b(\theta_{d-1}+j)(u^{i+1}_{j+1},u^i_j)$ is, modulo the image of
$\delta^{i-1}_T$, congruent to an element $(0,v^i_j)$, which of course
is also in the
kernel of $\delta^i_T$.
So it suffices to show that any such kernel element $(0,v^i_j)$ satifies
$b(\theta_{d-1}+j)(0,v^i_j)\in\im\delta^{i-1}_T$.

Since $\delta^i_T(0,v^i_j)=0$, we must have $\delta^iv^i_j=0$. Hence
$b(\theta_d+j)v^i_j=\delta^{i-1}a^{i-1}_j$ for some $a^{i-1}_j\in
\L^{i-1}_j$. Now 
$b(\theta_d+j)v^i_j=b(\theta_{d-1}+j)v^i_j+x_dQ\del_dv^i_j$
 for some $V_d$-homogeneous 
$Q\in F^0_H(D_n)\backslash F^{-1}_H(D_n)$. Therefore
$b(\theta_{d-1}+j)v^i_j=\delta^{i-1}a^{i-1}_j-x_dQ\del_dv^i_j$.

Since
$\delta^i(Q\del_dv^i_j)=-Q\del_d\delta^i(v^i_j)=0$, it follows that
$b(\theta_{d-1}+j)(0,v^i_j)=\del^{i-1}_T(Q\del_dv^i_j,a^{i-1}_j)$
and the induction claim is proved.

Now recall the inductive definition of
$\K^\bullet(\L^\bullet,x_1,\ldots,x_d)$, which together with the
induction claim shows that if the
cohomology of $V_d$-degree $j$ in $\L^\bullet$ is killed by
$b(\theta_d+j)$ then the cohomology of the complex
$\K^\bullet(\L^\bullet,x_1,\ldots,x_d)[j]$ is killed by $b^{(2^{d})}(j)$. 

The theorem now follows easily from the fact that $\C$ is a domain.\qed
\end{thm}

Theorem \ref{b-function-complex} 
leads to the computation of de Rham cohomology as follows.

Recall that we need to compute
$\Omega\otimes_{D_n}^LMV^\bullet$, and notice that 
$\Omega\otimes^L_{D_n}(-)$ is for any given complex quasi-isomorphic to
$\K^\bullet(-,\del_1,\ldots,\del_n)$. 
In order to cope with the problem that we'd like to compute
the complex $\K^\bullet(MV^\bullet,\del_1,\ldots,\del_n)$ and not
$\K^\bullet(MV^\bullet,x_1\ldots,x_n)$, we shall make use of the
Fourier transform. Namely, 
$H^i(\K^\bullet(MV^\bullet,\del_1,\ldots,\del_n))$ is isomorphic to 
$H^i(\K^\bullet(\tilde{MV}^\bullet,x_1,\ldots,x_n))$, where
$\tilde{MV}^\bullet$ is the image of the complex $MV^\bullet$ under
the Fourier transform $x\to\del, \del\to -x$.

Then $\tilde{MV}^\bullet$ may be replaced by a
free $V_n$-strict complex as $A^\bullet$ constructed in section
\ref{section-free-complex}. In particular, $A^i$ is then a
$\gr^\bullet_H(D_n)$-module and the cohomology of $A^\bullet$ is
holonomic (\cite{K2}) and therefore specializable to the origin. 

Let $b(s)\in \C[s]$ be the polynomial found in
subsection \ref{subsection-global-b-function} for $d=n$.
Then $b(\theta_d+j)\gr_H^j(H^iA^\bullet[\m_i])$ is zero
according to (\ref{b-function-kills}). Therefore by theorem
\ref{b-function-complex}, $b(j)$ kills the degree $j$ pieces of the
cohomology of $\K^\bullet(A^\bullet,x_1,\ldots,x_n)$. In other words,
if $k_0,k_1$ are integers with $b(s)=0$ only if $s\in [k_0,k_1]\cap
\Z$, then $\gr_H^j\K^\bullet(A^\bullet[\m_\bullet],x_1,\ldots,x_n)$ is
exact if 
$j\not\in[k_0,k_1]\cap\Z$.
Let $\tilde \Omega=F(\Omega)=D_n/(x_1,\ldots,x_n)\cdot D_n$.

We need to make a convention about the $V_n$-filtration on tensor
products with $\tilde\Omega$. 
If $A[\m]$ is a free $H$-graded $D_n$-module with shift
vector $\m$ then $\tilde\Omega\otimes_{D_n}A[\m]$ is filtered by
$F^i_H[\m](\tilde\Omega\otimes_{D_n}A):=\{\bar
P\otimes_{D_n}Q|V_ndeg(P)+V_ndeg[\m](Q)\le i\}$. Note that as $\tilde
\Omega$ equals $\C[\del_1,\ldots,\del_n]$ as right $D_n$-module,
$F^i_H[\m](\tilde\Omega\otimes_{D_n}A)$ equals
$\{(P_1,\ldots,P_{\dim_{D_n}A})|P_i\in\C[\del_1,\ldots,\del_n],
\deg_\del(P_i)\le \m(i)\,\,\forall i\}$. 
\begin{thm}
\mylabel{alg-de-rham}
The cohomology of $\Omega\otimes_{D_n}^LMV^\bullet$ can 
be computed as follows:
\begin{enumerate}
\item Compute $MV^\bullet$ as in \cite{W-1}, algorithm 5.1 as complex
of finitely generated $D_n$-modules.
\item Compute a $V_n$-strict free complex 
\[
\cdots\to A^{r-2}[\m_{r-2}]\to
A^{r-1}[\m_{r-1}]\to 0
\] 
quasi-isomorphic to $\tilde{MV}^\bullet$, $\tilde{MV}^\bullet$
denoting the image of $MV^\bullet$ under the Fourier automorphism.
\item Relative to the induced filtration on
$H^i(A^\bullet[\m_\bullet])$ compute the $b$-functions
$b_i(s)$ for the restriction of $H^i(A^\bullet[\m_\bullet])$ 
to the origin
$x_1=\ldots=x_n=0$. 
\item Let $b(s)$ be the least common multiple of all the
$b_i$. Find integers $k_0,k_1$ with $(b(k)=0, k\in\Z)\Rightarrow
(k_0\le k\le k_1)$. 
\item $\Omega\otimes_{D_n}^LMV^\bullet$ is
quasi-isomorphic to the complex 
\begin{eqnarray}
&\cdots\to\frac{F^{k_1}_H[\m_i](\tilde\Omega\otimes_{D_n}A^i)}
	{F^{k_0-1}_H[\m_i](\tilde\Omega\otimes_{D_n}A^i)}\to
\frac{F^{k_1}_H[\m_{i+1}](\tilde\Omega\otimes_{D_n}A^{i+1})}
	{F^{k_0-1}_H[\m_{i+1}](\tilde\Omega\otimes_{D_n}A^{i+1})}\to
\cdots
&
\label{restricted-complex}
\end{eqnarray}
shifted $n$ spots to the right.
\end{enumerate}
\proof
We already remarked that
$H^i(\K^\bullet(MV^\bullet,\del_1,\ldots,\del_n))$ is isomorphic to 
$H^i(\K^\bullet(\tilde{MV}^\bullet,x_1,\ldots,x_n))$, and so we 
only need to show that the latter can be computed from
(\ref{restricted-complex}). 

Almost tautologically, 
\[
\cdots\to \gr^k_H[\m_i]\tilde{MV}^i\to\cdots\to \gr^k_H[\m_{r-2}]\tilde{MV}^{r-2}
\to \gr^k_H[\m_{r-1}]\tilde{MV}^{r-1}\to 0
\]
is quasi-isomorphic to 
\[
\cdots\to \gr^k_H[\m_i]A^i\to\cdots\to \gr^k_H[\m_{r-2}]A^{r-2}
\to \gr^k_H[\m_{r-1}]A^{r-1}\to 0
\]
and therefore the same is true after application of
$\K^\bullet(-,x_1,\ldots,x_n)$. Moreover,
$\K^\bullet(\gr^\bullet_HA^i[\m_i],x_1,\ldots,x_n)[k]$ is
quasi-isomorphic to $\gr^k_H(\tilde \Omega\otimes_{D_n} A^i)[\m_i]$
shifted $n$ places to the right. Hence
$\K^\bullet(\gr^\bullet_H\tilde{MV}^\bullet[\m_\bullet],x_1,\ldots,x_n)[k]$
is 
quasi-isomorphic to 
\begin{eqnarray}
\label{complex}
\phantom{0000}&\cdots\to
\gr^k_H[\m_{r-2}](\tilde\Omega\otimes_{D_n} A^{r-2})
\to \gr^k_H[\m_{r-1}](\tilde\Omega\otimes_{D_n} A^{r-1})\to 0&
\end{eqnarray}
shifted $n$ places to the right. 
By theorem \ref{b-function-complex} this latter complex is exact for
all $k\not\in[k_0,k_1]\cap\Z$.

Observe that $F^k_H[0](\tilde\Omega\otimes_{D_n} D_n)$ is 
zero if $k<0$ for obvious
reasons. Therefore $F^k_H[\m_i](\tilde\Omega\otimes_{D_n} A^i)=0$ for
$k<\min_j\{\m_i(j)\}$. This together with
the exactness of (\ref{complex}) for $j<k_0$ forces 
$F^j_H[\m_i]H^i(\tilde\Omega\otimes_{D_n} A^\bullet)$ to vanish for
$j<k_0$. Similarly, $F^j_H[\m_i]H^i(\tilde\Omega\otimes_{D_n} A^\bullet)=
F^{j+1}_H[\m_i]H^i(\tilde\Omega\otimes_{D_n} A^\bullet)$ for $j\geq k_1$.

Thus the cohomology of $\tilde\Omega\otimes_{D_n} 
A^\bullet$ is captured by the quotient of  
$F^{k_1}[\m_i]H^i(\tilde\Omega\otimes_{D_n} 
A^\bullet)$
modulo
$F^{k_0-1}[\m_i]H^i(\tilde\Omega\otimes_{D_n} A^\bullet)$. 
The theorem follows.\qed
\end{thm}
\begin{rem}\label{this_remark}
\hfill
\par{\ref{this_remark}.1.}
The quotient $
\frac{F^{k_1}_H[\m_i](\tilde\Omega\otimes_{D_n}A^i)}
	{F^{k_0-1}_H[\m_i](\tilde\Omega\otimes_{D_n}A^i)}$ should be
thought of as a set of polynomials in $\del_1,\ldots,\del_n$ of
degrees bounded between $k_1-\m_i(j)$ and $k_0-\m_i(j)$.
\par{\ref{this_remark}.2.}
Since $\K^{r-i}(\gr^\bullet_H[\m_\bullet]A^\bullet,x_1,\ldots,x_n)$
involves only terms from $A^{r-1},\ldots,A^{r-i}$, the following  statement 
can be made: if $MV^\bullet$ is exact at $r-i$ and beyond, then
$\Omega\otimes_{D_n}^LMV^\bullet$ is exact at $j\geq r-i$. That follows
by considering $b(s)\cong 1$ which kills the last $i-1$ cohomology terms in
$MV^\bullet$, and inspecting the proof of theorem
\ref{b-function-complex} one sees that then $1$ also kills the
last $i-1$ cohomology terms in
$\K^\bullet(A^\bullet[\m_\bullet],x_1,\ldots,x_n)$.
\end{rem}
It follows the well-known 
\begin{cor}
$H^i_{dR}(U,\C)=0$ if $i\geq n+\cd(f_1,\ldots,f_r)$.
\end{cor}

\section{De Rham cohomology with support}
\mylabel{section-local}
Let $Y, Z$ be two Zariski-closed subsets of $X$. In this section we
are concerned with finding an algorithm that computes the de Rham
cohomology groups $H^\bullet_{dR,Z}(U,\C)$ 
of $U=X\backslash Y$ with coefficients in $\C$ and
support in $Z$.

$H^\bullet_{dR,Z}(U,\C)$ is defined as follows. Recall the de Rham
complex $\Omega^\bullet(U)$ on $U$. Usual de Rham cohomology is
defined as the hypercohomology of $\Omega^\bullet(U)$ and not
surprisingly de Rham cohomology with supports is defined as the
hypercohomology with supports in $Z$ 
of $\Omega^\bullet(U)$. In other words,
$H^\bullet_{dR,Z}(U,\C)=H^\bullet(R\Gamma_Z(U,\Omega^\bullet(U)))$. 

As was pointed out by Hartshorne, there is a natural exact sequence 
\begin{equation}
\label{long-local-seq}
\cdots\to H^i_{dR,Z}(U,\C)\to H^i_{dR}(U,\C)\to H^i_{dR}(U\backslash
Z,\C)
\to H^{i+1}_{dR,Z}(U,\C)\to\cdots
\end{equation}
which indicates that $H^\bullet_{dR,Z}(U,\C)$ measures the change in
cohomology due to the removal of $Z\cap U$ from $U$. 

For the entire section let us assume that $Y=\var(F), F=(f_1,\ldots,f_r)$ and
$Z=\var(G), G=(g_1,\ldots,g_s)$. Write $F\cdot G=\{f_i\cdot g_j\}$. 
As before we will write $F_I$ for
$\prod_{i\in I}f_i$ and so on. 
There is a natural map of Mayer-Vietoris
complexes 
\[
MV^\bullet(F\cdot G,F)\to MV^\bullet(F\cdot G)
\]
given by the natural projection 
\[
\bigoplus_{|I|+|J|+|K|=l}R_{F_I\cdot G_J}\otimes_RR_{F_K}\to
\bigoplus_{|I|+|J|=l} R_{F_I\cdot G_J}
\]
sending each summand with $|K|>0$ to zero. This map
corresponds to the embedding $X\backslash (Y\cup Z)\into
X\backslash Y$. Clearly the map is surjective and the kernel is the
subcomplex of $MV^\bullet(F\cdot G,F)$ consisting of those pieces
which contain at least one factor from $F$. It is not hard to check
that this kernel is exactly $MV^\bullet(F)\otimes_R C^\bullet(F\cdot
G)$, $C^\bullet(F\cdot G)$ being the \v Cech complex to $F\cdot G$
given by 
$\bigotimes_{i,j}(0\to R\stackrel{nat}{\longrightarrow}R_{f_i\cdot
g_j}\to 0)$.

Notice that the sequences 
\[
0\to (MV^\bullet(F)\otimes_R C^\bullet(F\cdot G))^i\to MV^i(F\cdot G,F)\to
MV^i(F\cdot G)\to 0
\]
are all split exact. Let $A^\bullet $ be a resolution of $\Omega$ as
right $D_n$-module, for example $A^\bullet$ could be the global sections
of the de Rham complex on $X$. Then there is a sequence of complexes 
\begin{equation}
\label{short-complex-seq}
A^\bullet\otimes_{D_n} (MV^\bullet(F)\otimes_R C^\bullet(F\cdot G))
\to A^\bullet \otimes_{D_n} MV^\bullet(F\cdot G,F)
\to A^\bullet \otimes_{D_n} MV^\bullet (F\cdot G)
\end{equation}
with split exact rows. In other words, we have a short exact sequence of
complexes. 

As was explained in previous sections, the cohomology of $A^\bullet
\otimes_{D_n} MV^\bullet(F\cdot G,F)$ is $H^i_{dR}(X\backslash Y,\C)$ while the
cohomology of $A^\bullet \otimes_{D_n} MV^\bullet(F\cdot G)$ is
$H^i_{dR}(X\backslash (Y\cup Z),\C)$ and the map on cohomology 
is induced by the
natural inclusion.

Comparison of the long exact sequence (\ref{long-local-seq}) with 
the long exact sequence that results from the short exact sequence of
complexes (\ref{short-complex-seq}) shows that the cohomology of
$A^\bullet\otimes_{D_n} (MV^\bullet(F)\otimes_R C^\bullet(F\cdot G))$ is
exactly $H^i_{dR,Z}(X\backslash Y,\C)$.

Computationally this is of course horrible: de Rham cohomology of
$X\backslash Y$ and $X\backslash Z$ comes from the Mayer-Vietoris
complex of $F$ and $G$ while here we have (approximately) the
Mayer-Vietoris complex of $F\cup F\cdot G$. We shall try to improve
this situation now. As a first step in that direction we point out
that the long exact sequence (\ref{long-local-seq}) shows that for
complements of affine closed varieties 
$H^i_{dR,Z}(X\backslash Y,\C)$ is in fact nothing but the
relative cohomology group $H^i(X\backslash Y, X\backslash (Y\cup
Z);\C)$. 

Consider the space $X\backslash (Y\cap Z)$ and its open covering by
the two sets $X\backslash Y$ and $X\backslash Z$. It follows from
\cite{G-H}, Example 17.1, that this is an exact triad for homology
with integer coefficients, and from \cite{E-S}, Theorem 11.4, that the
same holds for cohomology with coefficients in $\C$. This means that
the natural inclusion of pairs
\[
(X\backslash Y,X\backslash (Y\cup Z))\into (X\backslash (Y\cap
Z),X\backslash Z)
\]
induces an isomorphism between $H^i(X\backslash Y,X\backslash (Y\cup
Z);\C)$ and $H^i(X\backslash (Y\cap
Z),X\backslash Z;\C)$. This in turn implies that 
instead of
$H^i_{dR,Z}(X\backslash Y,\C)$ 
we may calculate
$H^i_{dR,Z}(X\backslash (Y\cap Z),\C)$ 
since the groups are isomorphic. 

Now consider the natural projection of complexes 
\[
MV^\bullet(F,G)\to MV^\bullet(G)
\]
given by
$\bigoplus_{|I|+|J|=l}R_{F_I}\otimes_RR_{G_J}\to\bigoplus_{|I|=l}R_{F_I}$
induced by the inclusion $X\backslash Z\into X\backslash (Y\cap
Z)$. As before, this induces a short exact sequence of complexes
\[
0\to MV^\bullet(F)\otimes_R C^\bullet(G)\to MV^\bullet(F,G)\to
MV^\bullet (G)\to 0.
\]
Tensoring over $D_n$ with the resolution $A^\bullet$ from above we discover that
the cohomology of $A^\bullet\otimes_{D_n} (MV^\bullet(F)\otimes_R
C^\bullet(G))$ is $H^\bullet_{dR,Z}(X\backslash (Y\cap Z))\cong
H^\bullet_{dR,Z}(X\backslash Y)$. Now the complexity is down to the
level of computing $H^\bullet_{dR}(X\backslash (Y\cap Z))$. It follows
\begin{alg}
\mylabel{de-rham-with-support}
Input: polynomials $F=\{f_1,\ldots,f_r\}$ defining $Y$ and
$G=\{g_1,\ldots,g_s\}$ defining $Z$; $i\in {\mathbb N}$.

Output: The de Rham cohomology groups of $U=X\backslash Y$ 
with support in $Z$,
$H^i_{dR,Z}(X\backslash Y,\C)$, which equal 
the relative cohomology groups
 $H^i(X\backslash Y, X\backslash(Y\cup
Z);\C)$.

Begin
\begin{enumerate}
\item Compute the complex $MV^\bullet(F)\otimes_R C^\bullet(G)$ as a
complex of left $D_n$-modules as in \cite{W-1}, algorithm 5.1.
\item Compute a free $V_n$-strict complex $A^\bullet$ quasi-isomorphic to
the image of $MV^\bullet(F)\otimes_R C^\bullet(G)$ under the Fourier
automorphism as in section \ref{restrict-complex}.
\item Find the cohomology groups of $A^\bullet$ and compute the
$b$-functions of the cohomology groups under the given shifts. Let
$k_0$ and $k_1$ be lower and upper bounds of the roots of the
$b$-functions. 
\item Replace each $D_n$ in $A^\bullet$ by
$k[\del_1,\ldots,\del_n]=\tilde\Omega$  and
restrict the complex to the components between $V_n$-degree $k_0-1$ and
$k_1$.
\item Take $n-i$-th cohomology of the resulting complex of
$\C$-vectorspaces and return it.
\end{enumerate}

End.
\end{alg}

\bibliography{bib}
\bibliographystyle{abbrv}

\end{document}